\documentclass[12pt]{amsart}
\usepackage{amsthm}
\usepackage{epsf}
\usepackage{amssymb}
\usepackage{latexsym}
\usepackage{amsmath}
\long\def\forget#1\forgotten{}
\forget
\newtheorem{theorem}{Theorem} 
\newtheorem{definition}[theorem]{Definition}

\newtheorem{lemma}[theorem]{Lemma}

\newtheorem{proposition}[theorem]{Proposition}
\newtheorem{remark}[theorem]{Remark}
\forgotten

\def\R{{\Bbb R}}

\def\C{{\Bbb C}}

\def\Z{{\Bbb Z}}

\def \ta{\tau}
\def \ta1{\tau_1}

\def \s{\sigma}

\def \A{{\mathcal A}}
\def \B{{\mathcal B}}

\def \zovera {
    \mathop{\lower 10pt \hbox{${\buildrel{\displaystyle\bar{z}} \over {\scriptstyle{(a)}}} $}}
    {\lower 4pt \hbox{${\scriptstyle{ij}}$}} 
} 

\newcommand\begintable[1][] {{}}

\long\def\forget#1\forgotten{}

\newif\ifXY 
\XYtrue     
\newif\ifbigmatrices
\bigmatricestrue 
 \ifXY
 \usepackage{xy}
 \fi
 \ifXY
 \xyoption{all}
 \fi

\usepackage{amsthm, epsf, amssymb,xy,latexsym, amsmath}

\newtheorem{theorem}{Theorem} 

\newtheorem{cor}[theorem]{Corollary}
\newtheorem{lemma}[theorem]{Lemma}
\newtheorem{proposition}[theorem]{Proposition}

\def \1{^{-1}}
\def \2{^{-2}}
\def \3{^{-3}}
\def \A{{\mathcal A}}
\def \B{{\mathcal B}}
\def \Cc{{\mathcal C}}
\newcommand\proj{{\mathbb P}^2}
\newcommand\comp[1]{\proj\backslash #1}
\newcommand\fg[1]{\pi_1(\comp{#1})}

\def\C{{\mathbb C}}
\def\R{{\mathbb R}}
\def\Z{{\mathbb Z}}
\def\t{{\tau}}
\def\q{{\kappa}}
\def\s{{\sigma}}
\def\l{{\lambda}}
\def\L{{\Lambda}}
\def\tgamma{{\tilde{\gamma}}}
\def\ttheta{{\tilde{\theta}}}
\def\tTheta{{\tilde{\Theta}}}
\def\tGamma{{\tilde{\Gamma}}}
\def\ogamma{{\overline{\gamma}}}

\newcommand\Ref[1]{~(\ref{#1})}
\newcommand\oio{{^{(i)}}}

\begin{document}

\renewcommand{\subjclassname}{%
       \textup{2000} Mathematics Subject Classification}
\date{\today}

\address{Meirav Amram, Mathematisches Institut, Bismarck Strasse 1 1/2, 
Erlangen, Germany}
\email{meirav@macs.biu.ac.il, amram@mi.uni-erlangen.de}

\address{Mina Teicher, Mathematics Department, Bar-Ilan University, Ramat-Gan, Israel}
\email{teicher@macs.biu.ac.il}
\address{A. Muhammed Uluda\u g, Mathematics Department, Galatasaray University,
Ortak\"oy, Istanbul, Turkey}
\email{uludag@gsunv.gsu.edu.tr}
 \title {Fundamental groups of some quadric-line arrangements}
 \author{ Meirav Amram \and Mina Teicher \and A. Muhammed Uludag}

\begin{abstract}
 In this paper we obtain  presentations of fundamental groups of the
 complements of three quadric-line arrangements in $\proj$. 
The first arrangement is a smooth quadric $Q$ with $n$ tangent lines to $Q$, 
and the second one is a quadric $Q$ with 
 $n$ lines passing through a point $p\notin Q$.  The last arrangement 
consists of a quadric $Q$ with $n$ lines passing through a point 
$p\in Q$.   
 \end{abstract} 
\maketitle
\section{\bf  Introduction.}
This is the first of a series of articles in which we shall study the 
fundamental groups of complements of some quadric-line arrangements. 
In contrast with the extensive literature on line arrangements and 
the fundamental groups of their complements, 
(see e.g. \cite{orliksalomon}, \cite{garberteicher} \cite{suciu}), 
only a little  
known about the quadric-line arrangements (see \cite{moishezonteicher5}).
The present article is dedicated to the computation of the fundamental 
groups of the complements of three infinite families of such arrangements.
A similar analysis for the quadric-line arrangements up to degree six 
will be done in our next paper.

\par
Let $C\subset \proj$ be a plane curve and $*\in \comp{C}$ a base point. 
By abuse of language we will call the group 
$\pi_1({\proj\backslash C},*)$ the \textit{fundamental group of C}, 
and we shall 
frequently omit base points and write $\fg{C}$.
One is interested in the group $\fg{C}$ mainly for the study of
the Galois coverings $X\rightarrow \proj$ branched along $C$.
Many interesting surfaces have been constructed as branched Galois coverings 
of the plane, for example
for the arrangement $\A_3$ in Figure 1 below, there are Galois coverings 
$X\rightarrow\proj$ branched along $\A_3$ such that 
$X\simeq {\mathbb P}^1\times {\mathbb P}^1$, or $X$ is an abelian surface, 
a K$3$ surface, or a quotient of the two-ball ${\mathbb B}_2$ 
(see \cite{yoshida}, \cite{holzapfel}, \cite{uludag1}). 
Moreover, some line arrangements 
defined by unitary reflection groups studied in \cite{orlik} 
are related to $\A_3$ via orbifold coverings. 
For example, if ${\mathcal L}$ is the line arrangement given by the equation 
$$
xyz(x+y+z)(x+y-z)(x-y+z)(x-y-z)=0
$$
then the image of ${\mathcal L}$ under the branched covering map 
$[x:y:z]\in \proj\rightarrow[x^2:y^2:z^2]\in\proj$ is the arrangement $\A_3$,
see \cite{uludag1} for details.
\par
The standard tool for fundamental group computations is the  
Zariski-van Kampen algorithm \cite{zariski}, \cite{vankampen},
see \cite{cheniot} for a modern approach.
We use a variation of this algorithm developed in~\cite{uludag3}
for computing the fundamental groups of 
real line arrangements and avoids lengthy monodromy computations.
The arrangements $\B_n$ and $\Cc_n$ discussed below are of fiber type, 
so presentations of their fundamental groups could be easily found as an 
extension of a free group by a free group. However, our 
approach has the advantage that it permits to capture the local 
fundamental groups around the singular points of these arrangements. 
The local fundamental groups are needed for the study of the singularities
of branched  of $\proj$ branched along these arrangements.
\par
In Section 2 below, we give fundamental group presentations 
and prove some immediate corollaries. 
In Section 3 we deal with the computations of fundamental
group presentations given in Section 2.

\section{\bf Results.}
Let $C\subset \proj$ be a plane curve and $B$ an irreducible component of $C$. 
Recall that a \textit{meridian} $\mu$ of $B$ in $\comp{C}$ 
with the base point $*\in \proj$ is a loop in $\comp{C}$ obtained by following 
a path $\omega$ with $\omega(0)=*$ and $\omega(1)$ belonging to a small 
neighborhood of a smooth point $p\in B\backslash C$, turning around $C$ in the 
positive sense along the boundary of a small disc $\Delta$ intersecting $B$ 
transversally at $p$, and then turning back to $*$ along $\omega$. 
The meridian $\mu$ represents a homotopy class in $\pi_1(\comp{C},*)$, 
which we also call a meridian of $B$. 
Any two meridians of $B$ in $\comp{C}$ are conjugate 
elements of $\fg{C}$ (see e.g.~\cite{lamotke}, 7.5), 
hence the meridians of irreducible components
of $C$ are supplementary invariants of the pair $(\proj,C)$. 
These meridians are specified in presentations of the fundamental group  
given below, they will be used in orbifold-fundamental 
group computations in \cite{uludag1}. 

\subsection{\bf The arrangement $\A_n$.}

\begin{figure}[h]
\begin{minipage}{\textwidth}
\begin{center}
\epsfbox{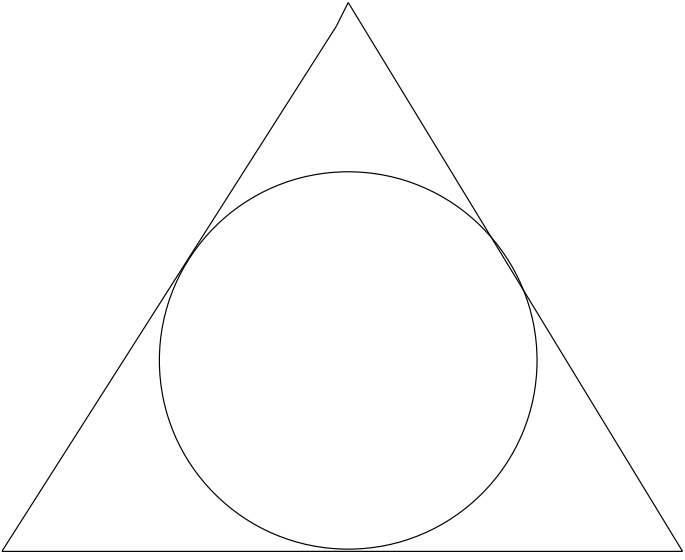}
\end{center}
\end{minipage}
\caption{The arrangement $\A_3$}\label{apoll}
\end{figure}
 

\begin{theorem}\label{an}
Let $\A_n:=Q\cup T_1\cup\dots\cup T_n$ be an arrangement consisting of a 
smooth quadric $Q$ with $n$ distinct tangent lines 
$T_1,\dots, T_n$. Then 
\begin{equation}\label{group1}
\fg{\A_n} \simeq \Biggl\langle 
\begin{array}{l}\t_1,\dots,\t_n,\\
 \q_1,\dots,\q_n 
\end{array}\; \left | 
\begin{array}{ll}
 \q_i = \t_i \q_{i-1} \t_i^{-1},\  2 \leq i \leq n \\
(\q_i \t_i)^2 = (\t_i \q_i)^2,\   1 \leq i \leq n \\
\left[\q_i^{-1} \t_i \q_i, \t_j \right] = 1,\  1 \leq i<j\leq n\\
\t_n \cdots \t_1 \q_1^2 = 1 
\end{array} 
\right . \Biggr\rangle 
\end{equation}
where $\q_i$ are meridians of $Q$ and $\t_i$ is a meridian of 
$T_i$ for $1\leq i\leq n$. Local fundamental groups around the singular points of $\A_n$ are generated by $<\q_i\1\t_i\q_i,\t_j>$ for the nodes $T_i\cap T_j$ 
and  by $<\q_i,\t_i>$ for the tangent points $T_i\cap Q$.
\end{theorem}
Part \textit{(i)} of the corollary below is almost trivial. 
Part \textit{(ii)} appears in~\cite{DOZ}, and part 
\textit{(iii)} was given in~\cite{degtyarev}. 
\begin{cor}
(i) One has: $\fg{\A_1}\simeq \Z$.\\
(ii) The group $\fg{\A_2}$ admits the presentation 
\begin{equation}\label{A2}
\fg{\A_2}\simeq \langle \t,\q\,|\, (\t\q)^2=(\q\t)^2 \rangle,
\end{equation}
where $\q$ is a meridian of $Q$ and $\t$ is a meridian of $T_1$. A meridian of
$T_2$ is given by $\q^{-2}\t\1$.\\
(iii) The group $\fg{\A_3}$ admits the presentation 
\begin{equation}
\fg{\A_3}\simeq \langle \t,\s,\q\,|\, 
(\t\q)^2=(\q\t)^2,\; (\s\q)^2=(\q\s)^2,\;[\s,\t]=1\rangle
\end{equation}
where $\s$, $\t$ are meridians of $T_1$ and $T_3$ respectively, 
and $\q$ is a meridian of $Q$. A meridian of $T_2$ is given by 
$(\q\t\q\s)\1$.
\end{cor}
A group $G$ is said to be \textit{big} if it contains a 
non-abelian free subgroup, and \textit{small} if $G$ is almost solvable. 
In  \cite{DOZ}, it was proved by V. Lin that the group\Ref{A2} is big. 
Below we give an alternative proof:
\begin{proposition}\label{big}
For $n>1$, the group $\fg{\A_n}$ is big.
\end{proposition}
\proof
A group with a big quotient is big. 
Since $\t_{n+1}$ is a meridian of $T_{n+1}$ in 
$\fg{\A_{n+1}}$, one has
$$
\fg{\A_n}\simeq \fg{\A_{n+1}}/\ll \t_{n+1}\gg,
$$ 
and it suffices to show that 
the group $\fg{\A_2}$ is big. In the presentation\Ref{A2}, 
applying the change of generators $\alpha:=\t\q$, $\beta:=\t$  gives
$$
\fg{\A_2}\simeq \langle \alpha, \beta\,|\, [\alpha^2,\beta]=1 \rangle.
$$
Adding the relations $\alpha^2=\beta^3=1$ to the 
latter presentation gives a surjection 
$\fg{\A_2}\twoheadrightarrow {\mathbb Z}/(2)*{\mathbb Z}/(3)$. 
Since the commutator subgroup 
of ${\mathbb Z}/(2)*{\mathbb Z}/(3)$ is the free group on two generators 
(see~\cite{dimca}), we get the desired result.

\subsection{The arrangement $\B_n$}\label{arrangementbn}
\begin{theorem}\label{bn}
Let $\B_n:=Q\cup T_1\cup T_2\cup L_1\cup\dots\cup L_n$ be an 
arrangement consisting of a smooth quadric $Q$ with $n+2$ distinct lines 
$T_1,T_2,L_1,\dots, L_n$ all passing through a point $p \notin Q$ 
such that $T_1$, $T_2$ are tangent to $Q$.
Then one has
\begin{equation}\label{group2}
\fg{\B_n} \simeq \Biggl\langle \t, \q,\l_1,\dots,\l_n \; \left |
\begin{array}{ll}
(\q \t)^2 = (\t \q)^2,\ \\
\left[\q,\l_i\right]=1,\ 1\leq i\leq n\\
\left[\t\1\q\t,\l_i\right]=1,\ 1\leq i\leq n
\end{array}
\right . \Biggr\rangle
\end{equation}
where $\t$ is a meridian of $T_1$, $\l_i$ is a meridian of $L_i$ for 
$1\leq i\leq n$, and $\q$ is a meridian of $Q$. A meridian $\s$ of $T_2$ 
is given by $\s:=(\l_n\dots \l_1\q^2\t)\1$.  
Local fundamental groups around the singular points of 
$\B_n$ are generated by $<\q, \l_i>$ and $<\t\1\q\t, \l_i>$ 
for the nodes $L_i\cap Q$,  by $<\q,\t>$ for the tangent point $T_1\cap Q$,
and by $<\q,\s>$ for the tangent point $T_2\cap Q$.
\end{theorem}

\begin{figure}[h]
\begin{minipage}{\textwidth}
\begin{center}
\epsfbox{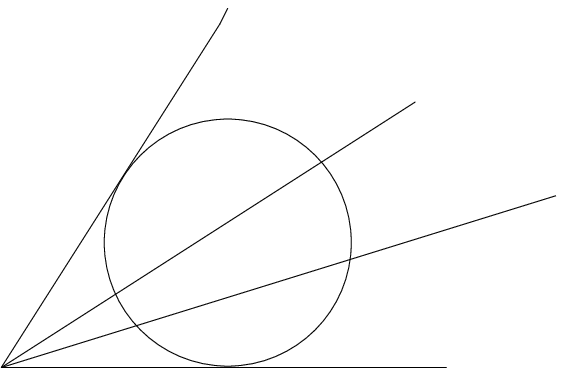}
\hspace{1cm}
\epsfbox {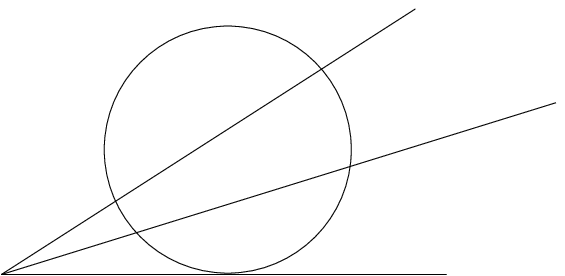}
\end{center}
\end{minipage}
\caption{Arrangements $\B_2$ and $B_2^\prime$}\label{b2}
\end{figure}

\begin{cor}\label{cor1} (i) Put $\B_n^\prime:=\B_n\backslash T_1$
and $\B_n^{\prime\prime}:=\B_n^\prime\backslash T_2$. Then 
\begin{equation}\label{bnprime}
\fg{\B_n^\prime}\simeq\fg{\B_{n+1}^{\prime\prime}}\simeq 
\langle \q,\l_1,\dots,\l_n\,|\, [\q,\l_i]=1,\,1\leq i\leq n \rangle
\end{equation}
\end{cor}
\noindent
\textit{Proof.} One has 
$\fg{\B_n^\prime}\simeq \fg{\B_n}/\ll \t\gg$. Setting $\t=1$ in  
presentation\Ref{group2} gives 
$$
\fg{\B_n^\prime}\simeq\langle \q, \l_1,\dots,\l_n,\,|\,
[\q,\l_i]=1\quad 1\leq i\leq n \rangle. 
$$
Setting $\t=1$ in the expression for a meridian $\s$ of $T_2$ given in 
Theorem~\ref{bn} shows that $(\l_n\dots \l_1\q^2)\1$ is a meridian of 
$T_2$ in $\fg{\B_n^\prime}$.
In order to find $\fg{\B_n^{\prime\prime}}$, it suffices to set 
$\l_n\dots \l_1\q^2=1$ in the presentation of $\fg{\B_n^\prime}$.
Eliminating $\l_n$ by this relation yields the presentation
$$
\fg{\B_n^{\prime\prime}}\simeq
\langle \q,\l_1,\dots,\l_{n-1}\,|\, [\l_i,\q]=[\l_{n-1}\dots\l_1\q^2,\q]=1
\rangle.
$$
Since the last relation above is redundant, we get the desired isomorphism
$\fg{\B_n^{\prime\prime}}\simeq \fg{\B_{n+1}^{\prime}}$. $\Box$

\medskip

Note that the groups $\fg{B_i^{\prime\prime}}$ are abelian for $i=0,1,2$.
Hence, the groups $\fg{B_i^\prime}$ are abelian for $i=0,1$. 
Otherwise, setting $\q=1$ in  presentation\Ref{bnprime} 
gives the free group  on $n-1$ generators, 
which shows that these groups are big. 
The groups $\fg{\B_n}$ are always big, since the arrangement $\B_0$ is 
the same as 
$\A_2$, and $\fg{\A_2}$ is big by Proposition~\ref{big}.
\subsection{The arrangement $\Cc_n$}\label{arrangementcn}
\begin{figure}[h]
\begin{minipage}{\textwidth}
\begin{center}
\epsfbox{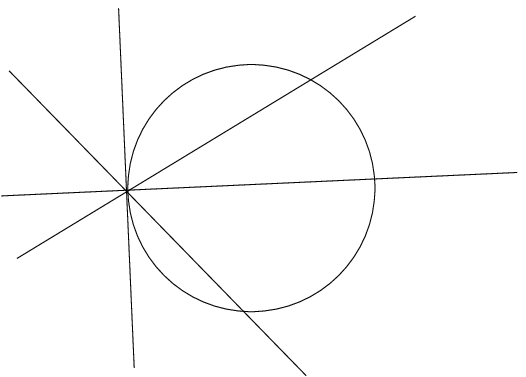}
\hspace{1cm}
\epsfbox {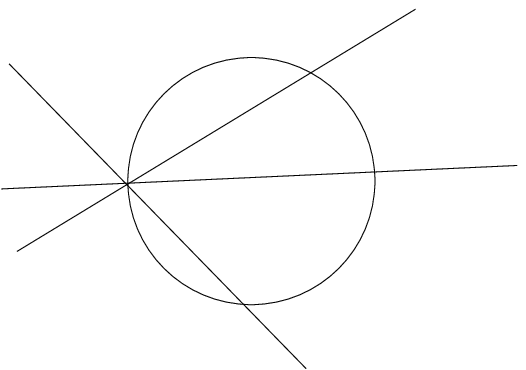}
\end{center}
\end{minipage}
\caption{Arrangements $\Cc_3$ and $\Cc_3^\prime$}\label{f1}
\end{figure}

\begin{theorem}\label{cn}
Let ${\Cc}_n:=Q\cup T\cup L_1\cup \dots\cup L_n$ 
be an arrangement consisting of a smooth quadric 
$Q$ with $n+1$ distinct lines $T,L_1\dots, L_n$, all passing 
through a point $p\in Q$ 
such that $T$ is tangent to $Q$. Then one has
\begin{equation}\label{group3}
\fg{\Cc_n}\simeq \langle \q,\l_1,\dots,\l_n\,|\, 
[\q,\l_i]=1\quad 1\leq i\leq n \rangle,
\end{equation}
where $\q$ is a meridian of $Q$ and $\l_i$ is a meridian of $L_i$ for
$1\leq i\leq n$. A meridian $\t$ of $T$ is given by
$\t:=(\l_n\dots\l_1\q^2)\1$. 
Local fundamental groups around the singular points of 
$\Cc_n$ are generated by $<\q, \l_i>$  
for the nodes $L_i\cap Q$, and by $<\t,\l_1,\dots,\l_n,\q>$ for the point $p$.
\end{theorem}

Note that the arrangement $\Cc_n$ is a degeneration (in the sense of Zariski)
of the arrangement $\B_n^\prime$ as the point $p$ approaches to $Q$. 
By Zariski's ``semicontinuity'' 
theorem of the fundamental group~\cite{zariski} (see also~\cite{dimca}), 
there is a surjection 
$\fg{\Cc_n}\twoheadrightarrow \fg{\B_n^\prime}$. 
In our case, this is also an injection: 
\begin{cor}\label{cor3}
(i) $\fg{\B_n^\prime}\simeq \fg{\Cc_n}$.\\
(ii) Put $\Cc_n^\prime:=\Cc_n\backslash T$.
Then $\fg{\Cc_n}\simeq \fg{\Cc_{n+1}^\prime}$.
\end{cor}
\noindent\textit{Proof.} Part (i) is obvious. The proof of 
part (ii) is same as the proof of Corollary~\ref{cor1}, (ii).
 
\section{\bf The arrangement $\A_n$}\label{anproof}
\begin{figure}[h]
\begin{minipage}{\textwidth}
\begin{center}
\epsfbox{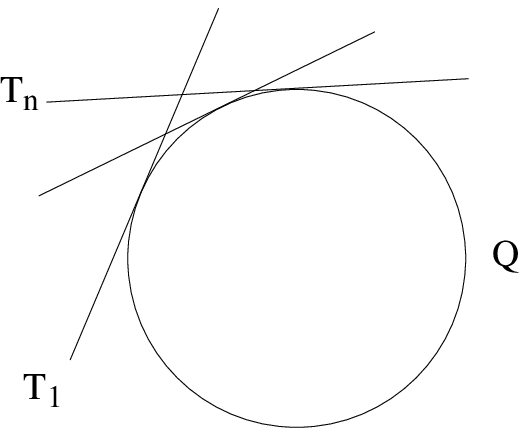}
\end{center}
\end{minipage}
\caption{}\label{}
\end{figure}

\medskip
It is easily seen that any two arrangements $\A_n$ 
with fixed $n$ are isotopic. 
In particular, the groups $\fg{\A_n}$ are isomorphic. 
Hence one can take as a model of the arrangements $\A_n$ the quadric 
$Q$ defined by 
$x^2+y^2=z^2$, where $[x:y:z]\in\proj$ is a 
fixed coordinate system in $\proj$. 
Pass to the affine coordinates in $\C^2\simeq \proj\backslash \{z=0\}$. 
Choose real numbers $x_1,\dots, x_n$ such that 
$-1<x_1<x_2\dots<x_n<0$, and define $y_i$ to be the positive solution of 
$x_i^2+y_i^2=1$ for $1\leq i\leq n$. 
Put $t_i:=(x_i,y_i)\in Q$, and take $T_i$ to be the tangent line 
to $Q$ at the point $t_i$ (see Figure 4).
\par
Let $pr_1:\C^2\backslash \A_n \rightarrow \C$ be the first projection. 
The base of this projection will be denoted by $B$. 
Put $F_x:=pr_1^{-1}(x)$, and denote by $S$ the set 
of singular fibers of $pr_1$. 
It is clear that if $F_x\in S$, then $x\in [-1,1]$. 
There are three types of singular fibers:\\
(i) The fibers $F_1$ and $F_{-1}$, corresponding to the `branch points' 
$(-1,0)$ and $(1,0)$. \\
(ii) The fibers $F_{x_i}$ $(1\leq i\leq n)$ corresponding to the `tangent points' 
$t_i=(x_i,y_i)= T_i\cap Q$.\\
(iii) The fibers $F_{a_{i,j}}$ $(1\leq i\neq j\leq n)$ corresponding to the nodes 
$n_{i,j}=(a_{i,j},b_{i,j}):=T_i\cap T_j$.
One can arrange the lines $T_i$ such that 
$$
-1<x_1<a_{1,2}<a_{1,3}<\dots<a_{1,n}<x_2<a_{2,3}<\dots <x_n<1
$$ 
\par
Identify the base $B$ of the projection $pr_1$ with the line 
$y=-2\subset \C^2$.
Let $N$ be the number of singular fibers and let 
$-1=s_1<s_2<\dots <s_{N-1}<s_N=1$ be the elements of $S\cap B$  
(so that $s_2=x_1$, $s_3=a_{1,2}$, $s_4=a_{1,3}$, and so on)
In $B$, take small discs $\Delta_i$ around the points $s_i$, and denote 
by $c_i$, $d_i$ ($c_i<d_i$) the points 
$\partial\Delta_i\cap \R$ for $1\leq i\leq N$ (see Figure 5).  

\begin{figure}[h]
\begin{minipage}{\textwidth}
\begin{center}
\epsfbox{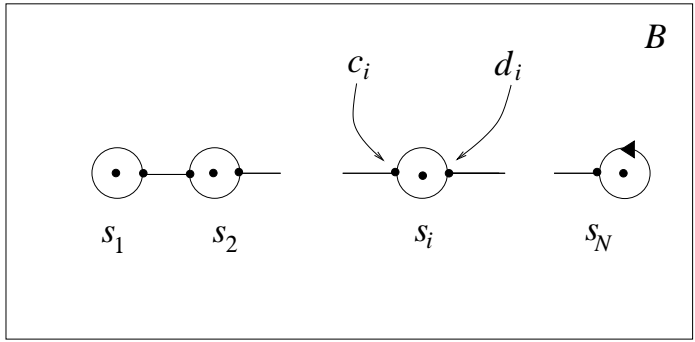}
\end{center}
\end{minipage}
\caption{The base $B$}\label{}
\end{figure}

\par
Put $B_1:=[c_1,c_2]\cup \Delta_1$ and for 
$2\leq i\leq N$ let $B_i:=[c_1,c_{i+1}]\cup \Delta_1\cup\dots\cup \Delta_i$.
Let $X_i:=pr^{-1}(B_i)$ be the restriction of the fibration $pr$ to $B_i$.
Let 
$$
A_i:=\Delta_i\cup \partial\bigl(\{\Im(z)\leq 0,\, c_2\leq \Re(z)\leq c_i\}\backslash (\Delta_2\cup \Delta_3\cup \dots \cup\Delta_{i-1})\bigr)
$$
and let $Y_i:=pr^{-1}(A_i)$ be the 
restriction of the fibration $pr$ to $A_i$. (see Figure 6).

\begin{figure}[h]
\begin{minipage}{\textwidth}
\begin{center}
\epsfbox{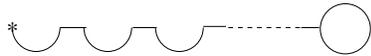}
\end{center}
\end{minipage}
\caption{The space $A_i$}\label{}
\end{figure}

Clearly, $X_{i}=X_{i-1}\cup Y_{i}$ for $2\leq i\leq N$.
We will use this fact to compute the groups 
$\pi_1(X_i,*)$ recursively, where $*:=(c_2,-2)$ is the base point. 
For details of the algorithm we apply below, see~\cite{uludag3}.
\par
 Identify the fibers of $pr_1$ with $F_0$ via the second projection 
$pr_2:(x,y)\in\C^2\rightarrow y\in\C$. 
 In \textit{each one} of the fibers $F_{c_i}$ (respectively $F_{d_i}$) 
take a basis for $\pi_1(F_{c_i}, -2)$ (respectively for $\pi_1(F_{d_i},-2)$) 
as in Figure 7 (for $F_{d_i}$, just replace $\gamma$'s by 
$\theta$'s in Figure 7). We shall denote these basis by the vectors 
 $\Gamma_i:=[\gamma_{1}\oio,\dots,\gamma_{n+2}\oio]$, 
 (respectively 
$\Theta_i:=[\theta_{1}\oio,\dots,\theta_{n+2}\oio]$).

\begin{figure}[h]
\begin{minipage}{\textwidth}
\begin{center}
\epsfbox{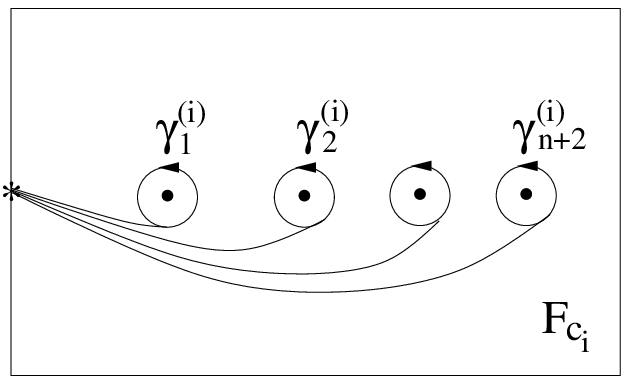}
\end{center}
\end{minipage}
\caption{}\label{}
\end{figure}

\par
Let $\nu_i\subset B_i\subset B$ be a path starting at $\nu_i(0)=c_2$,
ending at $\nu_i(1)=c_i$ and such that
$$
\nu_i([0,1])=\partial\bigl(\{\Im(z)\leq 0,\, c_2\leq \Re(z)\leq c_i\}\backslash (\Delta_2\cup \Delta_3\cup \dots \cup\Delta_{i-1})\bigr)
$$ 
Similarly, let $\eta_i\subset B_i\subset B$ be a path starting at 
$\eta(0)=c_2$, ending at $\eta(0)=d_i$ and such that 
$$
\eta_i([0,1])=\partial\bigl(\{\Im(z)\leq 0,\, c_2\leq \Re(z)\leq d_i\}\backslash (\Delta_2\cup \Delta_3\cup \dots \cup \Delta_{i})\bigr)
$$
For $2\leq i\leq N$ and $1\leq j\leq n+2$  each loop 
$\tgamma_j\oio:=\nu_i\cdot \gamma_j\oio\cdot \nu_i^{-1}$ 
represents a homotopy class in 
$\pi_1(X_{i},*)$, where $*:=(c_2,-2)$ is the base point. Similarly,
each loop $\ttheta_j\oio:=\eta_i\cdot\theta_i\cdot \eta_i^{-1}$ 
represents a homotopy class in $\pi_1(X_i,*)$. 
Denote $\tGamma_i:=[\tgamma_{1}\oio,\dots,\tgamma_{n+2}\oio]$, 
and $\tTheta_i:=[\ttheta_{1}\oio,\dots,\ttheta_{n+2}\oio]$.
\par
 It is well known that 
 the group $\pi_1(Y_i,*)$ has the presentation 
\begin{equation}
\label{pioneyi}
\langle \tgamma_{1}\oio,\dots, \tgamma_{n+2} \,|\,  
\tgamma_{j}\oio=M_i(\tgamma_{j}\oio),\, 1\leq j\leq n+2 \rangle
\end{equation}
 where $M_i:\pi_1(F_{c_i},-2)\rightarrow \pi_1(F_{c_i},-2)$ is the monodromy 
 operator around the singular fiber above $s_i$.  
 It is also well known that if it is the branches of $\A_n$ corresponding to  
 the loops $\tgamma_{k}\oio$ and $\tgamma_{k+1}\oio$ that meet above $s_i$,
then the only non-trivial relation in\Ref{pioneyi} is 
$\tgamma_{k}\oio=\tgamma_{k+1}\oio$ in case of a branch point,
$[\tgamma_{k}\oio,\tgamma_{k+1}\oio]=1$ in case of a node, and 
$(\tgamma_{k}\oio\tgamma_{k+1}\oio)^2=(\gamma_{k}\oio\gamma_{k+1}\oio)^2$
in case of a tangent point.  
\par
Now suppose that the group $\pi_1(X_{i-1},*)$ is known, with generators
$\tGamma_2$.
Recall that $X_i=X_{i-1}\cup Y_i$.
In order to find the group $\pi_1(X_i,*)$, one has to  express 
the base $\tTheta_i$ in terms of the base $\tGamma_i$. 
Adding to the presentation of $\pi_1(X_{i-1})$ 
the relation obtained by writing the relation of $\pi_1(Y_i)$ in the new base 
then yields a presentation of $\pi_1(X_i)$.
Note that, since the space $Y_i$ is eventually glued to $X_{i-1}$, 
it suffices to find an expression of $\tGamma_i$ in terms of the 
base $\tGamma_2$ \textit{in} the group $\pi_1(X_{i-1},*)$.

\par
Since all the points of $\A_n$ above the interval 
$[d_{i-1}, c_i]$ are smooth and real, one has

\medskip
\noindent
\textbf{Fact.} The loops $\ttheta_j^{(i-1)}$ and $\tgamma_j\oio$ are
homotopic in $X_i$ (or in $Y_i$) for $2\leq i\leq N$ and $1\leq j\leq n+2$. 
In other words, the bases $\tTheta_{i-1}$ and $\tGamma_i$ are homotopic.

\medskip
In order to express the base $\tTheta_i$ in terms of the base
$\tGamma_i$ the following lemma will be helpful. 
\begin{figure}[h]
\begin{minipage}{\textwidth}
\begin{center}
\epsfbox{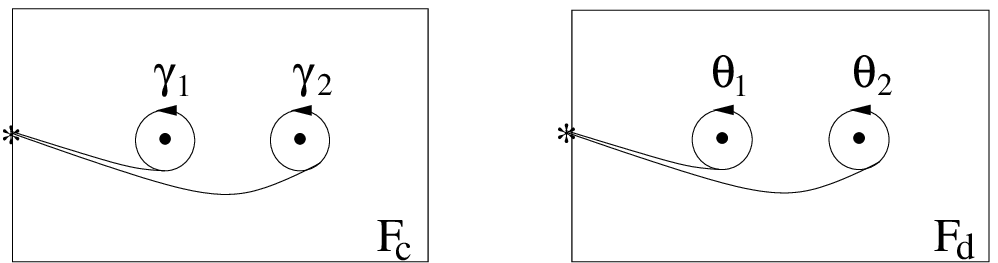}
\end{center}
\end{minipage}
\caption{}\label{}
\end{figure}

\begin{lemma}\label{otherside}
Let $C_k:x^2-y^{k+1}=0$  be an $A_k$ singularity, where $k=1$ or $k=3$.
Put $D:=\{(x,y)\,:\, |x|\leq 1, \,|y|\leq 1\}$ and let 
$pr_1:=(x,y)\in D\backslash C_k\rightarrow (x,-1)$ be the first projection. 
Denote by $F_x$ the fiber of $pr_1$ above $(x,-1)$. 
Identify the fibers of $pr_1$ via the second projection. 
Let $-1<c<0$ be a real number, and put $d:=-c$.
In $F_c$ (respectively in $F_d$) take a basis 
$\Gamma:=[\gamma_1,\gamma_2]$ for $\pi_1(F_c,-1)$ 
(respectively a basis $\Theta:=[\theta_1,\theta_2]$ for 
$\pi_1(F_d,-1)$) as in Figure 8.
Let $\eta$ be the path $\eta(t):=ce^{\pi it}$, and put 
$\ttheta_i:=\eta\cdot\theta\cdot\eta\1$ for $i=1,2$. Then
$\gamma_i$, $\ttheta_i$ are loops in $D\backslash C_k$ based at $*:=(c,-1)$,
and one has
\par
(i) If $k=1$, then $\ttheta_1$ is homotopic to $\gamma_2$, and 
$\ttheta_2$ is homotopic to $\gamma_1$, in other words,  
$$
\tTheta=[\gamma_2,\gamma_1].
$$
\par
(ii) If $k=3$, then 
$$
\tTheta=
[\gamma_2\gamma_1\gamma_2^{-1}, \gamma_1^{-1}\gamma_2 \gamma_1].
$$
\end{lemma}
\noindent\textit{Proof.}
Since $\pi_1(D\backslash C_2)$ is abelian, part (i) is obvious. 
For part (ii), note that the points of intersection 
$F_{\eta(t)}\cap C_4$ are $y_1:=c^2e^{2\pi it}$ and $y_2:=-c^2e^{2\pi it}$.
Hence, when we move the fiber $F_c$ over $F_d$ along the path $\eta$,
$y_1$ and $y_2$ make one complete turn around the origin 
in the positive sense. The loops $\gamma_1$, $\gamma_2$ are transformed 
to loops $\ogamma_1,{\ogamma}_2\subset F_d$
as in Figure 9. It follows that the loop
$\eta\cdot{\ogamma}_i\cdot\eta\1$ is homotopic to $\gamma_i$ 
for $i=1,2$. This homotopy can be constructed explicitly as follows:
Let $\Phi_{\eta(t)}:F_c\rightarrow F_{\eta(t)}$ be the corresponding
Leftschez homeomorphism (see~\cite{moishezonteicher1}). Then
$$
H(s,t):=\left \{
\begin{array}{ll}
\eta(3s),& 0\leq s\leq t/3\\
\Phi_{\eta(t)}(\gamma_i(3(s-t/3)/(3-2t))),& t/3\leq s\leq 1-t/3\\
\eta(3(1-s)),&1-t/3\leq s\leq 1
\end{array}\right .
$$
gives a homotopy between $\gamma_i$ and ${\ogamma}_i$.
\begin{figure}[h]
\begin{minipage}{\textwidth}
\begin{center}
\epsfbox{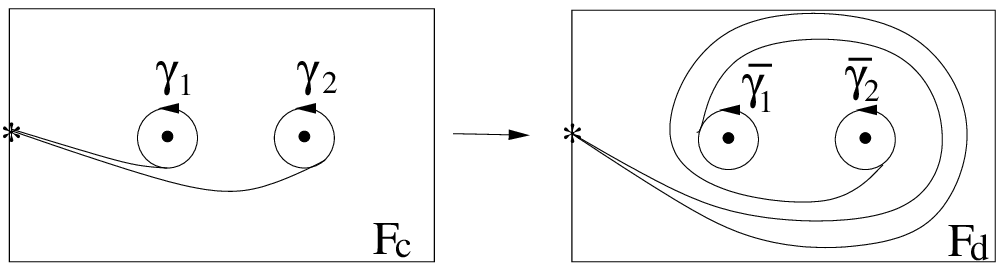}
\end{center}
\end{minipage}
\caption{}\label{}
\end{figure}

Expressing $\ttheta_i$ in terms of $\ogamma_i$, we get
$$
\ttheta_1=\ogamma_1\1\ogamma_2\1\ogamma_1\ogamma_2\ogamma_1=
\gamma_1\1\gamma_2\1\gamma_1\gamma_2\gamma_1,
$$
$$
\ttheta_2=\ogamma_1\1\ogamma_2\ogamma_1=\gamma_1\1\gamma_2\gamma_1.
$$
Since from the monodromy one has the relation  
$(\gamma_1\gamma_2)^2=(\gamma_2\gamma_1)^2$, the expression for  
$\ttheta_1$ can be simplified to get 
$\ttheta_1=\gamma_2\gamma_1\gamma_2\1$. $\Box$

\par 
Now we proceed with the computation of the groups $\pi_1(X_i)$. 
Clearly, the group $\pi_1(X_2)$ is generated by 
the base 
$$
\tGamma_2=[\gamma^{(2)}_1,\gamma^{(2)}_2,\dots ,\gamma^{(2)}_{n+2}]
$$
with the only relations 
\begin{equation}\label{rels1}
\gamma^{(2)}_1=\gamma^{(2)}_2
\end{equation}
 and 
\begin{equation}\label{rels2}
(\gamma^{(2)}_2\gamma^{(2)}_3)^2=(\gamma^{(2)}_3\gamma^{(2)}_2)^2.
\end{equation} 
Put 
$$
[\q_1,\q_1,\t_1,\dots,\t_n]:=\Gamma_2.
$$
Then relation\Ref{rels2} becomes
\begin{equation}
\quad (\q_1\t_1)^2=(\t_1\q_1)^2.
\end{equation} 
By Lemma~\ref{otherside} and the above Fact, one has
$$
\tGamma_3=\tTheta_2=
[\q_1, \t_1\q_1\t_1\1,\q_1\1\t_1\q_1,\t_2,\dots,\t_n].
$$
Since $s_3$ corresponds to the node $T_1\cap T_2$, the next relation is
\begin{equation}
\quad [\q_1\1\t_1\q_1,t_2]=1.
\end{equation}
Hence,
$$
\pi_1(X_3,*)\simeq
\langle \q_1,\t_1,\dots,\t_n\,|\, 
(\q_1\t_1)^2=(\t_1\q_1)^2,\quad [\q_1\1\t_1\q_1,\t_2]=1\rangle.
$$
By Lemma~\ref{otherside}, one has
$$
\tGamma_4=\tTheta_3=
[\q_1,\t_1\q_1\t_1\1,\t_2,\q_1\1\t_1\q_1,\t_3,\dots,\t_n].
$$
Since $s_4$ corresponds to the node $T_1\cap T_3$, one has the relation
$$
[\q_1\1\t_1\q_1,\t_3]=1.
$$
Hence,
$$
\pi_1(X_4,*)\simeq
\langle \q_1,\t_1,\dots,\t_n\,|\, 
(\q_1\t_1)^2=(\t_1\q_1)^2,\quad 
[\q_1\1\t_1\q_1,\t_2]=[\q_1\1\t_1\q_1,\t_3]=1\rangle.
$$
By Lemma~\ref{otherside}, one has 
$$
\tGamma_5=\tTheta_4=
[\q_1,\t_1\q_1\t_1\1,\t_2,\t_3,\q_1\1\t_1\q_1,\t_4,\dots,\t_n].
$$
Since $s_k$ corresponds to the node $T_1\cap T_{k-1}$ 
for $2\leq k\leq n+1$, repeating the above procedure 
gives the presentation
$$
\pi_1(X_{n+1},*)\simeq
\langle \q_1,\t_1,\dots,\t_n\,|\, 
(\q_1\t_1)^2=(\t_1\q_1)^2,\quad 
[\q_1\1\t_1\q_1,\t_k]=1\quad 2\leq k\leq n\rangle
$$
and
$$
\tGamma_{n+2}=\tTheta_{n+1}=
[\q_1,\q_2,\t_2,\t_3,\dots, \t_n, \q_1\1\t_1\q_1],
$$
where we put $\q_{i+1}:=\t_i \q_i \t_i\1$ for $1\leq i\leq n-1$.
\par
The next point $s_{n+2}$ corresponds to the tangent point 
$T_2 \cap Q$. This gives the relation 
\begin{equation}
(\q_2 \t_2)^2= (\t_2 \q_2)^2
\end{equation}
and 
$$
\tGamma_{n+2}=\tTheta_{n+1}=
[\q_1,\t_2\q_2\t_2\1,\q_2\1\t_2\q_2,\t_3,\dots, \t_n, \q_1\1\t_1\q_1].
$$ 
Now comes the $n-2$ points  $s_k$
corresponding to the nodes $T_2\cap T_{k-n}$ for $n+3\leq 2n+1$.
These give the relations
$$
[\q_2\1\t_2\q_2,\t_k]=1\quad 3\leq k\leq n.
$$
Hence, one has
$$
\pi_1(X_{n+1},*)\simeq
\langle \q_1,\q_2,\t_1,\dots,\t_n\,|\, 
$$
$$
\q_2=\t_1\q_1\t_1\1,\quad 
(\q_i\t_i)^2=(\t_i\q_i)^2,\quad
[\q_i\1\t_i\q_i,\t_k]=1\quad i< k\leq n,\quad i=1,2\rangle.
$$
We proceed in this manner until the last singular fiber $s_N$.
Since this is a branch point, the final relation is
\begin{equation}
\q_n=\q_1.
\end{equation}
This gives the presentation
\begin{equation}
\pi_1(X_N,*)\simeq \langle \t_i, \q_i, \; 1\leq i\leq n \left | 
\begin{array}{ll}
 \q_i = \t_i \q_{i-1} \t_i^{-1},\ 2 \leq i \leq n \\
(\q_i \t_i)^2 = (\t_i \q_i)^2,\  1 \leq i \leq n \\
\left[\q_i^{-1} \t_i \q_i, \t_j \right] = 1,\  1 \leq i<j\leq n\\
\q_1=\q_n 
\end{array} 
\right . \Biggr\rangle.
\end{equation}
Adding to this presentation of $\pi_1(X_N,*)$ 
the projective relation $\t_n\dots\t_1\q_1^2=1$ gives the presentation
\begin{equation}
\fg{\A_n} \simeq \langle \t_i, \q_i, \; 1\leq i\leq n \left | 
\begin{array}{ll}
 \q_i = \t_i \q_{i-1} \t_i^{-1},              &     2 \leq i \leq n \\
(\q_i \t_i)^2 = (\t_i \q_i)^2,                &     1 \leq i \leq n \\
\left[\q_i^{-1} \t_i \q_i , \t_j \right] = 1, &     1 \leq i<j\leq n\\
\t_n \cdots \t_1 \q_1^2 = 1,& \q_1=\q_n 
\end{array} 
\right . \Biggr\rangle 
\end{equation}
Note that the relation $\q_1=\q_n$ is redundant. 
Indeed, since $\q_i=\t_i\q_{i-1}\t_i\1$, one has
\begin{equation}\label{redundant}
\q_n=(\t_n\dots\t_1)\q_1(\t_n\dots\t_1)\1.
\end{equation}
But $\t_n\dots\t_1=\q^{-2}$ by the projective relation. Substituting this
in\Ref{redundant} yields the relation $\q_1=\q_n$. This finally
gives the presentation\Ref{group1} and proves Theorem \ref{an}.
Claims regarding the local fundamental groups around the singular points of 
$\A_n$ are direct consequences of the above algorithm.

\subsection{\bf Proof of Corollary 2.}

\smallskip\noindent
{\bf (i) The arrangement $\A_1$.} 
Writing down the presentation\Ref{group1} explicitly for $n=1$ gives
$$
\fg{\A_1}\simeq \Biggl\langle \q_1,\t_1 \left | \begin{array}{l}
(\q_1\t_1)^2 = (\q_1 \t_1)^2 \\
t_1 \q_1^2 = 1 \end{array} \right . \Biggr\rangle. 
$$
Eliminating $\t_1$ from the last relation shows that 
$\fg{A_1}\simeq\Z$.

\smallskip\noindent
{\bf (ii) The arrangement $\A_2$.}
Writing down the presentation\Ref{group1} explicitly for $n=2$ gives
$$
\fg{\A_2} \simeq \Biggl\langle \q_1,\q_2,\t_1, \t_2 \left |
\begin{array}{l}
(1) \ \q_2 = \t_1 \q_1 \t^{-1}_1 \\
(2) \ (\q_1 \t_1)^2 = (\t_1 \q_1)^2 \\
(3) \ (\q_2 \t_2)^2 = (\t_2 \q_2)^2 \\
(4) \ [\q_1^{-1} t_1 \q_1 , t_2] = 1\\
(5) \ \t_2 \t_1 \q_1^2 = 1 \end{array} \right . \Biggr\rangle.
$$
Eliminating $\q_2$ by (1) and $\t_2$ by (5) 
one easily shows that the relations 
(3) and (4) are redundant. 
This leaves (2) and gives the desired presentation.

\smallskip\noindent
{\bf (iii) The arrangement $\A_3$.}
Writing down the presentation\Ref{group1} explicitly for $n=3$ gives
$$
\fg{\A_3} \simeq 
\Biggl\langle 
\begin{array}{l}\q_1,\q_2,\q_3,\\ \t_1,\t_2, \t_3\end{array} 
\left | 
\begin{array}{ll}
(1) \ \q_2 = \t_1\q_1\t_1\1 &
(2) \ \q_3 = (\t_2\t_1)\q_1(\t_2\t_1)\1 \\
(3) \ (\q_1\t_1)^2 = (\t_1 \q_1)^2 &
(4) \ (\q_2\t_2)^2 = (\t_2 \q_2)^2 \\
(5) \ (\q_3\t_3)^2 = (\t_3 \q_3)^2 &
(6) \ [\q_1^{-1} \t_1 \q_1,\t_2] = 1 \\
(7) \ [\q_1\1 \t_1 \q_1,\t_3] = 1 &
(8) \ [\q_2\1 \t_2 \q_2, \t_3 ] = 1 \\
(9)\ \t_3 \t_2 \t_1 \q_1^2 = 1
\end{array} \right . \Biggr\rangle.
$$
Eliminate $\q_2$ by (1), $\q_3$ by (2), and $\t_2$ by (9).
It can be shown that the relations (4), (6) and (8) are consequences of the 
remaining relations. The relation (5) becomes 
$(\q_1\t_3)^2=(\t_3\q_1)^2$. This gives the presentation
$$
\fg{\A_3}\simeq 
\langle \q_1,\t_1,\t_3\,|\, (\q_1\t_1)^2=(\t_1\q_1)^2,
\quad (\q_1\t_3)^2=(\t_3\q_1)^2,\quad [\q_1\1\t_1\q_1,\t_3]=1\rangle.
$$
Finally, put $\q:=\q_1$, $\t:=\q_1\1\t_1\q_1$ and $\s:=\t_3$.
Then $\t_1=\q\t\q\1$, and the first relation in the above presentation becomes
$(\q^2\t\q\1)^2=(\q\t)^2\Rightarrow (\q\t)^2=(\t\q)^2$.
This gives the desired presentation.

\section{\bf The arrangement $\B_n$}
\begin{figure}[h]
\begin{minipage}{\textwidth}
\begin{center}
\epsfbox{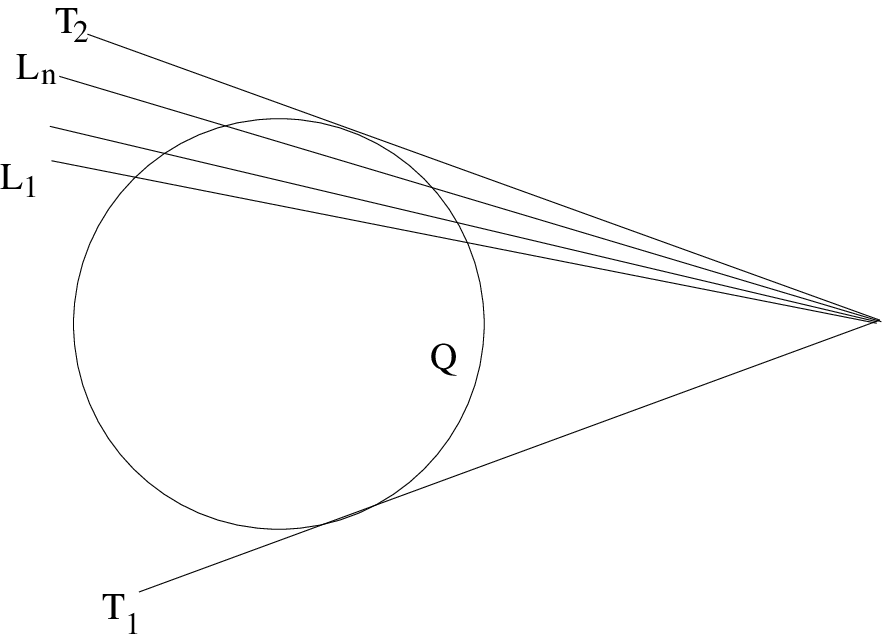}
\end{center}
\end{minipage}
\caption{}\label{}
\end{figure}

As in the case of the arrangements $\A_n$, it is readily seen that 
arrangements $\B_n$ are all isotopic to each other for fixed $n$, 
so one can compute $\fg{\B_n}$ from the following model for 
$\B_n$'s (see Figure 10): The quadric $Q$ is given by 
the equation $x^2+y^2=1$, and $p$ is the point $(2,0)$.
The lines $L_i$ intersect $Q$ above the $x-$ axis.
\par
The projection to the $x-$ axis 
has four types of singular fibers:\\
(i) The fibers $F_{1}$ and $F_{-1}$, corresponding to `branch points',\\
(ii) The fiber $F_{x}$ corresponding to the `tangent points' 
$(x,y)=t_1:=T_1\cap Q$ and $(x,-y)=t_2:=T_2\cap Q$,\\  
(iii) The fibers 
$F_{a_1},\dots, F_{a_n}$ ($-1<a_n<\dots a_1<x$) 
corresponding to the nodes $L_i\cap Q$ lying on the right of the tangent 
points and the fibers 
$F_{b_1},\dots, F_{b_n}$ ($x<b_1<\dots b_n<1$) corresponding to the 
nodes $L_i\cap Q$ lying on the left of the tangent points.\\
(iv) The fiber $F_{2}$, corresponding to the point $p$.\\
\par
In order to find the group $\fg{\B_n}$, we shall apply the same procedure
as in the computation of $\fg{\A_n}$.
Let $y\in \R$ be such that $-1<y<a_n$, and take $F_y$ to be the base fiber.
Let $s_1:=-1$, $s_{i+1}:=b_i$ for $1\leq i\leq n$, $s_{n+2}:=x$, 
$s_{n+2+i}:=a_{n+1-i}$ for $1\leq i\leq n$, and $s_{2n+3}:=1$, and
$s_{2n+4}=2$.
Take a basis 
$$
\tGamma_2:=[\t,\q_1,\q_2,\l_1,\dots,\l_n,\s]
$$
for $F_y$ as in Figure 7.  
\par
Since $s_1$ corresponds to a branch point, one has the relation
$\q_1=\q_2$. Put $\q:=\q_1=\q_2$. 
The point $s_2$ is a node, and yields the relation 
$[\q,\l_1]=1$, and one has
$$
\tGamma_3=[\t,\q,\l_1,\q,\l_2,\dots,\l_n,\s].
$$
Repeating this for the nodes $s_3,\dots, s_{n+1}$ gives the relations
$[\q,\l_i]=1$ for $1\leq i\leq n$, and one has
$$
\tGamma_{n+2}=[\t,\q, \l_1,\dots,\l_n,\q,\s]
$$
The monodromy around the fiber $F_{x}$ gives the relations
$(\t\q)^2=(\q\t)^2$ and $(\s\q)^2=(\q\s)^2$. One has
$$
\tGamma_{n+3}=
[\q\t\q\1,\t\1\q\t,\l_1,\dots,\l_n,\s\q\s\1,\q\1\s\q].
$$
Since the points $s_{n+3},\dots, s_{2n+2}$ corresponds to nodes,
one has the relations $[\l_i, \s\q\s\1]=1$, and
$$
\tGamma_{2n+3}=
[\q\t\q\1,\t\1\q\t,\s\q\s\1,\l_1,\dots,\l_n,\q\1\s\q].
$$
The branch point corresponding to $s_{2n+3}$ yields the relation
$$
\t\1\q\t=\s\q\s\1
$$ 
Together with the projective relation 
$\s\l_n\dots\l_1\q^2\t=1$ these relations already gives a presentation 
of $\fg{\B_n}$, since one can always ignore one of the singular fibers
when computing the monodromy (see \cite{uludag3}). 
\par
We obtained the presentation 
$$
\fg{\B_n} \simeq \Biggl\langle 
\L,\t, \q,\l_1,\dots,\l_n, \s \left |
\begin{array}{l}
(1) \ (\q \t)^2 = (\t \q)^2 \\
(2) \ (\q \s)^2 = (\s \q)^2 \\
(3) \ \t\1\q\t=\s\q\s\1\\
(4) \ [\q,\l_i] = 1 \quad 1\leq i\leq n\\
(5) \ [\s\q\s\1,\l_i]=1\quad 1\leq i\leq n\\ 
(6) \ \s\l_n\dots\l_1\q^2\t=1  \end{array} \right . \Biggr\rangle.
$$
Put $\L:=\l_n\dots\l_1$. 
Eliminating $\s$ by (7), it is easily seen that (3) is redundant.
Relation (2) becomes
$$
(\L\q^2\t\q\1)^2=(\q\1\L\q^2\t)^2\Rightarrow [\t\1\q\t,\L]=1.
$$
But this relation is a consequence of (4), so that (2) is also redundant.
Since $\t\1\q\t=\s\q\s\1$ by (3), the relation (5) can be written as
$[\t\1\q\t,\l_i]=1$. 
This gives the presentation\Ref{group2} and proves Theorem \ref{bn}.

\section{\bf The arrangement $\Cc_n$}
\begin{figure}[h]
\begin{minipage}{\textwidth}
\begin{center}
\epsfbox{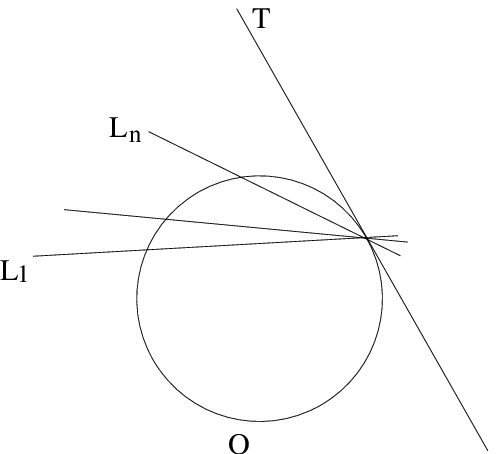}
\end{center}
\end{minipage}
\caption{}\label{}
\end{figure}

In order to compute the group, consider the model of $\Cc_n$ shown 
in Figure 11, where $Q$ is given by $x^2+y^2=1$.
Suppose that the second points of intersection of the 
lines $L_i$ with $Q$ lie above the $x-$ axis. 
As in the previous cases, take an initial base
$$
\tGamma_2:=[\q_1,\q_2,\l_1,\dots,\l_n,\t].
$$
The relation induced by the branch point is $\q_1=\q_2=:\q$.
The nodes of $\Cc_n$ will give the relations $[\q,\l_i]=1$ for 
$1\leq i\leq n$, and one has
$$
\tGamma_{n+2}=[\q,\l_1,\dots,\l_n,\q,\t].
$$
One can simplify the computation of 
 the monodromy around the complicated singular 
fiber as follows: Put $\Lambda:=\l_n\dots\l_1$.
By the projective relation one has
$\t\Lambda\q^2=1\Rightarrow \t=\q\2\Lambda\1$.
Hence, $[\q,\t]=1$. Since we also have $[\q,\l_i]=1$ for $1\leq i\leq n$, 
this means that when computing the monodromy
around this fiber, one can ignore the branch $Q$. This leaves
$n+1$ branches intersecting transversally, and the induced relation is 
(see~\cite{uludag3}) 
\begin{equation}\label{red}
[\t\Lambda,\lambda_i]=[\t\Lambda,\t]=1,
\end{equation}
and one has
$$
\tGamma_{n+3}=[\q,\q,....].
$$
The last relation induced by the branch point yields the trivial relation
$\q=\q$, as expected.
\par
Eliminating $\t$ shows that the relations\Ref{red} are redundant 
and gives the presentation
$$
\fg{C_n}\simeq
\langle \q,\l_1,\dots,\l_n\,|\,[\l_i,\q]=1\rangle.
$$
$\Box$

\bigskip\noindent
\textbf{Acknowledgements.} 
This work was partially supported by the Emmy Noether Research Institute for
 Mathematics  (center of the Minerva Foundation of Germany), the
Excellency Center
``Group  Theoretic Methods in the Study of Algebraic Varieties'' of the
Israel Science
 Foundation, and EAGER (EU network, HPRN-CT-2009-00099).

\end{document}